\documentclass[11pt]{amsart}

\usepackage{amsfonts,amssymb,amscd,amsmath,latexsym,amsbsy,amsthm}

\usepackage{amssymb}
\usepackage{amsfonts}
\usepackage{latexsym}

\theoremstyle{plain}
\newtheorem{theorem}{Theorem}

\newtheorem{lemma}[theorem]{Lemma}

\newtheorem{proposition}[theorem]{Proposition}

\newtheorem{corollary}[theorem]{Corollary}

\newtheorem{question}[theorem]{Question}

\newtheorem{conjecture}[theorem]{Conjecture}

\theoremstyle{definition}

\theoremstyle{remark}
\newtheorem*{remark}{Remark}

\numberwithin{theorem}{section}
\numberwithin{equation}{section}

\newcommand{\g}{\mathfrak{g}}

\newcommand{\wh}{\widehat{\mathfrak{h}}}

\newcommand{\ben}{\begin{enumerate}}
\newcommand{\een}{\end{enumerate}}

\newcommand{\QQ}{{\mathbb{Q}}}

\newcommand{\ZZ}{{\mathbb{Z}}}

\def\sym{{\mathfrak S}}

\newcommand{\solu}[1]{\begin{sol}{\bf (\ref{#1})}}

\begin{document}

\title{Groups and Lie algebras corresponding to the Yang-Baxter equations}

\author[L.B.]{Laurent Bartholdi}
\address{Institut de Math\'ematiques B (IMB), \'Ecole Polytechnique
  F\'ed\'erale de Lausanne (EPFL), CH-1015 Lausanne, Switzerland}
\email[L.B.]{laurent.bartholdi@epfl.ch}

\author[B.E.]{Benjamin Enriquez}
\address[B.E.]{IRMA, CNRS et Universit\'e de Strasbourg,
7 rue Ren\'e Descartes, F-67084 Strasbourg, France }
\email[B.E.]{enriquez@math.u-strasbg.fr}

\author[P.E.]{Pavel Etingof}
\address[P.E.]{Massachusetts Institute of Technology, Department of
  Mathematics 2-176, 77 Massachusetts Avenue, Cambridge, MA
  02139-4307, USA}
\email[P.E.]{etingof@math.mit.edu}

\author[E.R.]{Eric Rains}
\address[E.R.]{Department of Mathematics, One Shields Ave., University
  of California, Davis, CA 95616, USA}
\email[E.R.]{rains@math.ucdavis.edu}

\maketitle

\begin{abstract}
For a positive integer $n$ we introduce 
quadratic Lie algebras ${\mathfrak{tr}}_n$
${\mathfrak{qtr}}_n$ and finitely generated groups 
${\mathbf {Tr}}_n$, ${\mathbf {QTr}}_n$ naturally associated with the
classical and quantum Yang-Baxter equation, respectively. 

We prove that the universal enveloping algebras 
of the Lie algebras ${\mathfrak{tr}}_n$, 
${\mathfrak{qtr}}_n$ are Koszul, and compute 
their Hilbert series. We also compute the cohomology rings of these Lie
algebras (which by Koszulity are the quadratic duals of the 
enveloping algebras). 

We construct cell complexes which are 
classifying spaces for the groups ${\mathbf {Tr}}_n$
and ${\mathbf {QTr}}_n$, and show that the 
boundary maps in them are zero, which allows us to compute 
the integral cohomology of these groups. 

We show that the Lie algebras ${\mathfrak{tr}}_n$,
${\mathfrak{qtr}}_n$ map onto the associated graded algebras
of the Malcev Lie algebras of the groups ${\mathbf{Tr}}_n$,
${\mathbf{QTr}}_n$, respectively, and conjecture that 
this map is actually an isomorphism. 
(This conjecture was recently proved by P. Lee, \cite{L}). 
At the same time, we show that the groups ${\mathbf {Tr}}_n$
and ${\mathbf {QTr}}_n$ are not formal for $n\ge 4$. 
\end{abstract}
 
\section{Introduction}

In this paper we consider certain discrete groups and Lie
algebras associated to the Yang-Baxter equations.

Namely, we define the $n$-th quasitriangular 
Lie algebra ${\mathfrak{qtr}}_n$ to be generated by 
$r_{ij}$, $1\le i\ne j\le n$, with defining relations
given by the classical Yang-Baxter equation
\begin{equation}\label{cyb1}
[r_{ij},r_{ik}]+[r_{ij},r_{jk}]+[r_{ik},r_{jk}]=0
\end{equation}
for distinct $i,j,k$, 
and $[r_{ij},r_{kl}]=0$ for distinct $i,j,k,l$. 
We define the $n$-th triangular Lie algebra
${\mathfrak{tr}}_n$ by the same generators and relations, with the 
additional relations 
$$
r_{ij}=-r_{ji}.
$$ 
We define the $n$-th quasitriangular group $\mathbf{QTr}_n$ to be generated by 
$R_{ij}$, $1\le i\ne j\le n$, with defining relations
given by the quantum Yang-Baxter equation
\begin{equation}\label{cyb2}
R_{ij}R_{ik}R_{jk}=R_{jk}R_{ik}R_{ij}.
\end{equation}
and $R_{ij}R_{kl}=R_{kl}R_{ij}$ if $i,j,k,l$ are distinct. 
We define the $n$-th triangular group 
$\mathbf{Tr}_n$ by the same generators and relations, with the 
additional relations 
$$
R_{ij}=R_{ji}^{-1}.
$$ 
These groups and Lie algebras are trivial for $n=0,1$. 

\begin{remark} The groups $\mathbf{QTr}_n$
and $\mathbf{Tr}_n$ are also called pure virtual braid groups 
and pure flat braid groups, respectively, see \cite{L} and references 
therein.
\end{remark}

These definitions are motivated by the theory of 
quantum groups, as explained in Section 2. 

The main results of the paper are as follows. 

1. We prove that the universal enveloping algebras 
of the Lie algebras ${\mathfrak{tr}}_n$, 
${\mathfrak{qtr}}_n$ are Koszul, and compute 
their Hilbert series. We also compute the cohomology rings of these Lie
algebras (which by Koszulity are the quadratic duals of the 
enveloping algebras). 

2. We construct classifying spaces of the groups ${\mathbf {Tr}}_n$
and ${\mathbf {QTr}}_n$. More specifically, a classifying space  
for the group ${\mathbf{Tr}}_n$ can be obtained by gluing
faces of the $(n-1)$-th permutohedron corresponding to the same
set partition, and a similar construction works for ${\mathbf
{QTr}}_n$. Moreover, the boundary maps in the resulting cell complexes 
are both zero, which allows one to compute the cohomology of the groups 
${\mathbf{Tr}}_n$, ${\mathbf{QTr}}_n$ with integer
coefficients. 

3. We show that the Lie algebras ${\mathfrak{tr}}_n$,
${\mathfrak{qtr}}_n$ map onto the associated graded algebras
of the Malcev Lie algebras of the groups ${\mathbf{Tr}}_n$,
${\mathbf{QTr}}_n$, respectively. 

4. The quantum group intuition suggests a conjecture 
that these maps are isomorphisms; in other words, that the ranks of
the lower central series quotients for the groups ${\mathbf{Tr}}_n$,
${\mathbf{QTr}}_n$ are equal to the dimensions 
of the homogeneous components of the Lie algebras 
${\mathfrak{tr}}_n$, ${\mathfrak{qtr}}_n$.\footnote{We used the computer system
``GAP'' \cite{GAP4} to show that the conjecture was
true for ${\mathbf {Tr}}_4$ up to degree 7.}

We note, however, that the groups ${\mathbf {Tr}}_n$ for $n\ge 4$
are not formal in the sense of Sullivan~\cite{Su} (i.e.\ their Malcev Lie
algebras are not isomorphic to their rational holonomy Lie algebras). 
The same holds for the groups ${\mathbf {QTr}}_n$, $n\ge 4$.  
Summarizing, we may say that the properties 
of the groups ${\mathbf {Tr}}_n$, ${\mathbf {QTr}}_n$ 
are similar to those of the pure cactus group $\Gamma_n$ 
studied in \cite{EHKR}. 

{\bf Remark.} 1. We warn the reader that the published version of this paper, \cite{BEER},
contains a serious error. Namely, in \cite{BEER},  
the above conjecture for ${\mathbf{Tr}}_n$  
is stated as a theorem (Theorem 2.3), and a proof is given, which is 
incorrect. 
Namely, the proof rests on Proposition 5.1 of \cite{BEER}, 
which is false (we are grateful to  
P. Lee for discovering this error). As a result, the proof of Theorem 8.5 
given in \cite{BEER} contains a gap, as it rests on 
the incorrectly proved Theorem 2.3. 
Similarly, the proof of Proposition 6.1 of \cite{BEER} contains a gap, as it rests on the wrong 
Proposition 5.1. In the present version, these errors are corrected: Propositions 5.1 
and 6.1 are deleted, and Theorems 2.3 and 8.5 are stated as conjectures 
(Conjectures \ref{isomo} and \ref{cohring} below, 
respectively).

2. Luckily, Conjectures \ref{isomo} 
and \ref{cohring} were recently proved by P. Lee, \cite{L},
which effectively corrects the errors in \cite{BEER}. In fact, 
he also proved Conjectures \ref{isomo1} 
and \ref{cohring1}. 

\section{Lie algebras and groups corresponding to the Yang-Baxter equations}

\subsection{$I$-objects}

Let $I$ be the category of finite sets where morphisms are inclusions. 
Similarly, let $J$ be the category of ordered finite sets, where morphisms are
increasing inclusions. 

An $I$-object (respectively, a $J$-object) 
of a category ${\mathcal C}$ is a covariant functor
$I\to {\mathcal C}$ (respectively, $J\to {\mathcal C}$). 
Thus, an $I$-object (respectively, a $J$-object) is the same thing 
as a sequence of objects $X_1,X_2,\dots $ in $\mathcal C$ 
and a collection of maps $X(f): X_m\to X_n$ for every injective
(respectively, strictly increasing) 
map $f: [m]\to [n]$, such that $X(f)X(g)=X(fg)$.
\footnote{Here $[n]$ denotes the set $\lbrace{1,\dots ,n\rbrace}$.}
A morphism between $I$- and between $J$-objects is, 
by definition, a morphism of functors. 

Obviously, $J$ is a subcategory of
$I$ (with the same isomorphism classes of objects but fewer 
morphisms). Thus every $I$-object is also a $J$-object. 

{\bf Example.} Let $A$ be a unital associative algebra. Then we can define 
the $I$-algebra $T(A)$, such that $T(A)_n=A^{\otimes n}$, and 
for any $a\in A$ and $1\le k\le m$, $T(A)(f)(a_k)=a_{f(k)}$, where 
$a_k$ denotes the element $1\otimes\dots\otimes
a\otimes\dots\otimes 
1\in A^{\otimes m}$,
with $a$ being in the $k$-th component. 

\subsection{$I$-Lie algebras associated to the classical Yang-Baxter 
equation}

Let us define two $I$-Lie algebras over $\QQ$, 
${\mathfrak{tr}}$ and ${\mathfrak{qtr}}$ 
(the triangular and quasitriangular Lie algebra). 
Namely, the Lie algebras ${\mathfrak{tr}}_n$ and ${\mathfrak{qtr}}_n$
have been defined above. Now 
for each injective map $f: [m]\to [n]$ we have the corresponding map 
${\mathfrak{(q)tr}}(f): {\mathfrak{(q)tr}}_m\to {\mathfrak{(q)tr}}_n$ given
by $f(r_{ij})=r_{f(i)f(j)}$, which gives ${\mathfrak{qtr}}$ and ${\mathfrak{tr}}$ 
the structure of $I$-Lie algebras. 

We can also define the corresponding universal enveloping
$I$-algebras 
$U({\mathfrak{qtr}})$, $U({\mathfrak{tr}})$
in the obvious way. 

This definition is motivated by the following proposition, 
whose proof is straightforward. 

\begin{proposition} Let $A$ be a unital associative algebra.
Then ${\rm Mor}_I({\mathfrak{qtr}},T(A))$ (in the category of $I$-Lie
algebras) is 
the set of elements $r\in A^{\otimes 2}$
satisfying the classical Yang-Baxter equation
$$
[r_{12},r_{13}]+[r_{12},r_{23}]+[r_{13},r_{23}]=0,
$$
where $r_{ij}$ is the image in $A^{\otimes 3}$ of $r$ through
the map $1\mapsto i,2\mapsto j$. Similarly, ${\rm
  Mor}_I({\mathfrak{tr}},T(A))$ is the set of skew-symmetric
elements $r\in A^{\otimes 2}$
satisfying the classical Yang-Baxter equation.
\end{proposition}

We have natural homomorphisms ${\mathfrak{tr}}_n\to 
{\mathfrak{qtr}}_n\to {\mathfrak{tr}}_n$, whose composition 
is the identity. Namely, the second (surjective) map is defined by sending
$r_{ij}$ to $r_{ij}$ for all $i,j$, while the first (injective) map 
is defined by the same condition, but only for $i<j$. Thus 
 ${\mathfrak{tr}}_n$ is a split quotient of  ${\mathfrak{qtr}}_n$. 
We also note that the injection ${\mathfrak{tr}}_n\to
{\mathfrak{qtr}}_n$ induces a map 
$\phi: {\rm Mor}_I({\mathfrak{qtr}},T(A))\to 
{\rm Mor}_J({\mathfrak{tr}},T(A))$, which 
is actually an isomorphism; in particular,
${\rm Mor}_J({\mathfrak{tr}},
T(A))$ is the set of elements $r\in
 A^{\otimes 2}$ (not necessarily skew-symmetric) 
satisfying the classical Yang-Baxter equation.

\subsection{$I$-groups associated to the quantum Yang-Baxter equation.} 

We can also define $I$-groups $\mathbf{QTr}$, $\mathbf{Tr}$, 
which are quantum analogs of the Lie algebras 
${\mathfrak{qtr}}$, ${\mathfrak{tr}}$. Namely, 
for each injective map $f: [m]\to [n]$ we have the corresponding map 
$\mathbf{(Q)Tr}(f): \mathbf{(Q)Tr}_m\to \mathbf{(Q)Tr}_n$ given by 
$f(R_{ij})=R_{f(i)f(j)}$, which gives $\mathbf{QTr}$ and $\mathbf{Tr}$ 
the structure of $I$-groups. 

This definition is motivated by the following proposition, 
whose proof is straightforward. 

\begin{proposition} Let $A$ be a unital associative algebra.
Then ${\rm Mor}_I(\mathbf{QTr},T(A))$ (in the category of $I$-monoids) 
is the set of invertible elements $R\in A^{\otimes 2}$
satisfying the quantum Yang-Baxter equation
$$
R_{12}R_{13}R_{23}=R_{23}R_{13}R_{12},
$$
and ${\rm Mor}_I(\mathbf{Tr},T(A))$ is the set of elements $R\in A^{\otimes 2}$
satisfying the quantum Yang-Baxter equation and the unitarity
condition $R^{21}=R^{-1}$.
\end{proposition}

We have natural homomorphisms ${\mathbf{Tr}}_n\to 
{\mathbf{QTr}}_n\to {\mathbf{Tr}}_n$, whose composition 
is the identity. Namely, the second (surjective) map is defined by sending
$R_{ij}$ to $R_{ij}$ for all $i,j$, while the first (injective) map 
is defined by the same condition, but only for $i<j$. Thus 
 ${\mathbf{Tr}}_n$ is a split quotient of  ${\mathbf{QTr}}_n$. 
We also note that the injection ${\mathbf{Tr}}_n\to
{\mathbf{QTr}}_n$ induces a map 
$\phi: {\rm Mor}_I({\mathbf{QTr}},T(A))\to 
{\rm Mor}_J({\mathbf{Tr}},T(A))$, which 
is actually an isomorphism; in particular,
${\rm Mor}_J({\mathbf{Tr}},
T(A))$ is the set of elements $R\in  A^{\otimes 2}$ 
satisfying the quantum Yang-Baxter equation (but not necessarily
the unitarity condition).

{\bf Remark.} Recall that if $E,F$ are sets, then a (partially
defined) function $f:E \to F$ is a pair $(D_f,\bar f)$ 
where $D_f \subset E$, and $\bar f:D_f \to F$  is a map. 
Define $I'$ (resp., $J'$) as the opposite of the category 
whose objects are finite (resp., ordered
finite) sets and morphisms are partially defined functions (resp., non-decreasing 
functions). Then $\mathfrak{qtr}$, $\mathfrak{tr}$
are $I'$-Lie algebras, and $\mathbf{QTr}$, $\mathbf{Tr}$
are $J'$-groups. A partially defined function $f:[n]\to [m]$ gives rise to
a morphism $\mathfrak{(q)tr}_m \to \mathfrak{(q)tr}_n$
by $r_{ij}\mapsto \sum_{i'\in f^{-1}(i), j'\in f^{-1}(j)}
r_{i'j'}$. Similarly, a non-decreasing 
partially defined function $f : [n]\to [m]$ gives rise to the
morphism $\mathbf{(Q)Tr}_m \to \mathbf{(Q)Tr}_n$, by 
$R_{ij} \mapsto \prod_{i'\in
f^{-1}(i),j'\in f^{-1}(j)}R_{i'j'}$ 
(where the product is taken in increasing order of $i',j'$).

\subsection{Prounipotent completions} 
For a discrete group $G$, let ${\rm Lie}(G)$ denote the Lie algebra 
of the $\QQ$-prounipotent completion of $G$ (i.e., the Malcev
Lie algebra of $G$), and let ${\rm grLie}(G)$ 
be the associated graded of this Lie algebra with respect to 
the lower central series filtration. If $G$ is an $I$-group,
these are $I$-Lie algebras. 

We have natural surjective homomorphisms of $I$-Lie algebras 
$\phi_{{\mathfrak{(q)tr}}}: {\mathfrak{(q)tr}}\to \operatorname{grLie}(\mathbf{(Q)Tr})$, 
given by the formula $r_{ij}\to \log R_{ij}$. 

\begin{conjecture} \label{isomo} (stated as 
Theorem 2.3 in the published version; 
now a theorem of P. Lee, \cite{L})  
The homomorphism $\phi_{{\mathfrak{tr}}}$ is an isomorphism. 
\end{conjecture}

\begin{conjecture}\label{isomo1} (now a theorem of P. Lee, \cite{L})
The homomorphism $\phi_{{\mathfrak{qtr}}}$ is an 
isomorphism. 
\end{conjecture}

Note that Conjecture \ref{isomo} follows from Conjecture \ref{isomo1},
since the group $\mathbf{Tr}_n$ is a split quotient 
of the group $\mathbf{QTr}_n$ (via $R_{ij}\mapsto R_{ij}\mapsto R_{ij}$),
and similarly for the corresponding Lie algebras. 

Conjecture \ref{isomo1} is motivated by the following theorem. 

\begin{theorem}
Let $K$ be the kernel of $\phi_{{\mathfrak{qtr}}}$. Then 
$K$ is annihilated by every morphism in ${\rm
Mor}({\mathfrak{qtr}},T(A))$ 
for any algebra $A$. 
\end{theorem}

\begin{proof} By a result of \cite{EK,EK:2}, any 
solution $r\in A\otimes A$ of the classical Yang-Baxter equation 
can be quantized to a solution 
$R=1+\hbar r+O(\hbar^2)$ of the quantum Yang-Baxter equation. 
This implies that any morphism ${\mathfrak{(q)tr}}\to T(A)$ can be
deformed to a morphism 
$\mathbf{(Q)Tr}\to T(A)$, which implies the required statement. 
\end{proof}

{\bf Remark.} 
It is obvious that ${\mathfrak{tr}}_3$ is a free product 
of abelian Lie algebras 
$\QQ^2 * \QQ$, and ${\mathbf{Tr}}_3$ is a free product
$\ZZ^2 * \ZZ$; therefore ${\rm Lie}
({\mathbf{Tr}}_3)={\mathfrak{tr}_3}$. 
However, we have checked using ``Magma'' that 
${\rm Lie}(\mathbf{Tr}_n)$ is not isomorphic to ${\mathfrak{tr}_n}$ for
$n=4$, already modulo elements of degree $5$. 
Since we have split injections ${\mathfrak{tr}_n}\to {\mathfrak{tr}_{n+1}}$
and ${\mathbf{Tr}_n}\to {\mathbf{Tr}_{n+1}}$, the same statement
holds for $n>4$. This implies that 
the group $\mathbf{Tr}_n$ is not formal for $n\ge 4$ (i.e. its
classifying space is not a formal topological space, see \cite{Su}).
Since we have a split injection ${\mathbf {Tr}}_n\to {\mathbf
{Qtr}}_n$, the same holds\footnote{For basics about formality of groups, 
see e.g. \cite{PS}.} for the groups ${\mathbf {QTr}_n}$,
$n\ge 4$. 

\section{The Koszulity and Hilbert series of
$U({\mathfrak{tr}_n})$, $U({\mathfrak{qtr}_n})$.} 

One of the main results of this paper is the following theorem. 

\begin{theorem}\label{main}
(i) The algebras $U({\mathfrak{tr}}_n)$, $U({\mathfrak{qtr}}_n)$ are Koszul. 

(ii) The Hilbert series of these algebras are equal to 
$1/P_{{\mathfrak{(q)tr}}_n}(-t)$, where $P_{{\mathfrak{(q)tr}}_n}(t)$ are the 
polynomials with the following exponential generating functions:

\begin{equation}\label{Tr}
1+\sum_{n=1}^\infty P_{{\mathfrak{tr}}_n}(t)\frac{u^n}{n!}=e^{\frac{e^{tu}-1}{t}},
\end{equation}

\begin{equation}\label{QTr}
1+\sum_{n=1}^\infty P_{{\mathfrak{qtr}}_n}(t)\frac{u^n}{n!}=e^{\frac{u}{1-tu}}.
\end{equation}

(iii) The polynomials $P_{{\mathfrak{(q)tr}}_n}(t)$ are given by the following 
explicit formulas: 
$$
P_{{\mathfrak{tr}}_n}(t)=\sum_{k=1}^n\frac{1}{k!}
\Big( \sum_{i=0}^{k-1}(-1)^i\binom{k}{i}(k-i)^n \Big) t^{n-k}.
$$
(the palindrome of Bell's exponential polynomial, \cite{Be});
$$
P_{{\mathfrak{qtr}}_n}(t)=\sum_{p=0}^{n-1}\binom{n-1}{p}\frac{n!}{(n-p)!}t^p.
$$
\end{theorem}

The proof of this theorem is given in the next section.

\section{Proof of Theorem \ref{main}}

Part (iii) of the theorem follows by direct computation from 
part (ii), so we prove only parts (i) and (ii). 

\subsection{The triangular case}

Recall that if $B$ is a quadratic algebra, then the quadratic dual
$B^!$ is the quadratic algebra with generators $B[1]^*$ 
and relations $B[2]^\perp \subset B[1]^*\otimes B[1]^*$. 

Consider the quadratic dual algebra $A_n$ to
$U({\mathfrak{tr}}_n)$. Denote by $a_{ij}$ the set of generators of $A_n$
dual to the generators $r_{ij}$ of $U({\mathfrak{tr}}_n)$ (so 
$a_{ij}$ are defined for distinct $i,j\in [n]$, and $a_{ij}=-a_{ji}$).

\begin{lemma}\label{dualpres} The algebra $A_n$ is 
the supercommutative algebra generated by 
odd generators $a_{ij}$, $1\leq i\neq j\leq n$, with defining relations 
$a_{ij} + a_{ji} = 0$, and  
$$
a_{ij}a_{jk}=a_{jk}a_{ki}
$$
for any three distinct indices $i,j,k$. 
\end{lemma} 

\begin{proof} Let $\omega=\sum a_{ij}r_{ij}\in A_n\otimes {\mathfrak{tr}}_n$. The relations 
of $A_n$ can be written as the Maurer-Cartan equation $[\omega,\omega]=0$ (where the
commutator is taken in the supersense). 
Taking components of this equation, 
we get the relations for $a_{ij}$. 
\end{proof} 

Lemma \ref{dualpres} allows us to easily find a basis of $A_n$. 
Namely, define a monomial in $A_n$ to be reduced if
it is of the form $a_{i_1i_2}a_{i_1i_3}\dots a_{i_1i_m}$, 
with $i_1<i_2<\dots <i_m$. The support of this monomial is 
the set $\lbrace{i_1,\dots ,i_m\rbrace}$, and the root label is $i_1$. 

\begin{proposition}\label{basis} Products of reduced monomials 
with disjoint supports (in 
the order of increasing the root labels) form a basis in $A_n$.
\end{proposition}

\begin{proof}
Take any monomial in $A_n$. If it is not a product of 
reduced ones with disjoint supports, then 
it has a quadratic factor of the form $a_{ij}a_{jk}$, where $j>i$ or $j>k$. 
Using the relations, we can then replace it with another 
quadratic monomial, so that the total sum of labels is reduced. 
This implies that products of reduced monomials span $A_n$. 
The fact that these products are linearly 
independent is easy, since all relations are binomial. 
\end{proof}

\begin{corollary}\label{koz}
(i) The elements $a_{ij}a_{jk}-a_{jk}a_{ki}$ with $k<i,j$ form a
quadratic Gr\"obner basis for supercommutative algebras
\footnote{About Gr\"obner bases for supercommutative algebras,
see e.g. \cite{MV}.} 
for $A_n$ (for any ordering of monomials in which the sum of
labels is monotonically nondecreasing).

(ii) $A_n$ is Koszul. 

(iii) The Hilbert polynomial of $A_n$ is $P_{{\mathfrak{tr}}_n}(t)$. 
\end{corollary}

\begin{proof}
(i) follows from the Proposition \ref{basis}. (ii) follows from (i) 
since any supercommutative algebra with a quadratic Gr\"obner basis is Koszul 
(see e.g.\ \cite{Yu}, Theorem 6.16).
To prove (iii), note that Proposition \ref{basis} implies that $\dim A_n[k]$ 
is the number of partitions of the set $[n]$ into $n-k$ nonempty
parts (see \cite{wilf:gf}), so the result follows from 
standard combinatorics. 
\end{proof}

Corollary \ref{koz} and the standard theory of Koszul algebras imply 
Theorem \ref{main} in the triangular case. Indeed, the dual 
of a Koszul algebra is Koszul, and the Hilbert series of a Koszul 
algebra and its dual are related by the equation $p(t)q(-t)=1$. 

\subsection{The quasitriangular case}

The proof in the quasitriangular case is analogous although a bit more complicated.
Let us split 
$r_{ij}\in {\mathfrak{qtr}}_n$ into a symmetric and skew-symmetric part: 
$r_{ij}=t_{ij}+\rho_{ij}$, where $t_{ij}$ is symmetric and $\rho_{ij}$ 
is skewsymmetric in $i,j$. Then the defining relations for ${\mathfrak{qtr}}_n$ take the form:
$$
[t_{ij},t_{ik}+t_{jk}]=0,\ [t_{ij},\rho_{ik}+\rho_{jk}]=0,
$$
$$
[\rho_{ij},\rho_{jk}]+
[\rho_{jk},\rho_{ki}]+[\rho_{ki},\rho_{ij}]=[t_{ij},t_{jk}],
$$
for distinct $i,j,k$, and 
$$
[\rho_{ij},\rho_{kl}]=[\rho_{ij},t_{kl}]=[t_{ij},t_{kl}]=0
$$
if $i,j,k,l$ are distinct. 

As before, consider the quadratic dual algebra $QA_n$ to $U({\mathfrak{qtr}}_n)$.
Denote by $a_{ij}$, $b_{ij}$ the set of generators of $QA_n$
dual to the generators $\rho_{ij}$ and $t_{ij}$ of $U({\mathfrak{qtr}}_n)$ (so 
$a_{ij}$, $b_{ij}$ are defined for distinct $i,j\in [n]$, and $a_{ij}=-a_{ji}$, $b_{ij}=b_{ji}$).

\begin{lemma}\label{dualpres1} The algebra $QA_n$ is 
the supercommutative algebra generated by 
odd generators $a_{ij},b_{ij}$ with defining relations 
$$
a_{jk}a_{ij}=a_{ki}a_{jk}=b_{ij}b_{jk}+b_{jk}b_{ki}+b_{ki}b_{ij},
$$
$$
a_{ij}b_{jk}=a_{ik}b_{jk}
$$
for any three distinct indices $i,j,k$, and 
$$
a_{ij}b_{ij}=0,
$$
for $i\ne j$. 
\end{lemma} 

\begin{proof} 
The proof is analogous to the proof of Lemma \ref{dualpres}. 
\end{proof} 

The algebra $QA_n$ has a filtration in which $\deg(b_{ij})=1$, 
and $\deg(a_{ij})=0$. In the algebra ${\rm gr}QA_n$, 
the graded versions of the above relations are satisfied. 
These graded versions are the same as the original relations, except for the first 
set of relations in Lemma \ref{dualpres1}, which is replaced by 
$$
a_{ij}a_{jk}=a_{jk}a_{ki},\ b_{ij}b_{jk}+b_{jk}b_{ki}+b_{ki}b_{ij}=0.
$$

Let $QA_n^0$ be the algebra with generators $a_{ij}$, $b_{ij}$,
whose defining relations are the graded version of the relations of $QA_n$.
We have a surjective homomorphism 
$QA_n^0\to {\rm gr}QA_n$ (we will show below that 
it is an isomorphism). Note that we have a split injection $OS_n\to QA_n^0$ from 
the Orlik-Solomon algebra of the braid arrangement (generated by $b_{ij}$)
to $QA_n^0$. 

The exterior algebra in $a_{ij}$ and $b_{ij}$ is graded, as a space,
by 2-step set partitions: partitions of $[n]$ into nonempty sets
$S_1,\dots ,S_l$, and then of each $S_p$ into nonempty subsets
$S_{pq}$, $q=1,\dots ,m_p$.  Namely, if we are given a monomial $M$ in
$a_{ij}$ and $b_{ij}$, we connect $i,j$ by a black edge if $b_{ij}$ is
present in $M$, and by a red edge if $b_{ij}$ is not present, but
$a_{ij}$ is present. Then we define the $S_p$ to be the connected
components of the obtained graph, and the $S_{pq}$ to be the connected
components of the graph of black edges only.

It is easy to see that the relations of $QA_n^0$ are compatible with
this grading, and thus that $QA_n^0$ also has a grading by 2-step set
partitions.  This fact allows us to find a basis of $QA_n^0$.

Namely, let $S=(S_{pq})$ be a 2-step set partition of $[n]$, and 
let $i_{pq}$ be the minimal element of $S_{pq}$. Let $i_p$ be the minimum of $i_{pq}$ over all $q$,
and $q_p$ be such that $i_{pq_p}=i_p$.  

For $T\subset [n]$, let $OS(T)$ be the Orlik-Solomon algebra generated 
by $b_{ij}$, $i,j\in T$. Let $\lbrace{b(T,s)\rbrace}$, $1\le s\le (|T|-1)!$, be the broken circuit basis 
of the top component of this algebra (see \cite{Yu}). 

Let $QA_n^0(S)$ be the degree $S$ part of $QA_n^0$.

\begin{proposition}\label{basis1} 
  The elements $\prod_{p=1}^r
  (\prod_{q=1}^{s_p}b(S_{pq},s_{pq})\prod_{q\ne q_p}a_{i_pi_{pq}})$
  for all $s_{pq}$ form a basis of $QA_n^0(S)$.
\end{proposition}

\begin{proof}
  The proof is analogous to the proof of Proposition \ref{basis}.
  Namely, it is easy to show that any monomial in $QA_n^0(S)$ can be
  reduced, using the relations, to a monomial from Proposition
  \ref{basis1}. On the other hand, it is clear that the monomials in
  Proposition \ref{basis} are linearly independent (this follows from
  compatibility of the relations with the grading by 2-step set
  partitions).
\end{proof}

\begin{corollary}\label{koz1}
(i) The algebra $QA_n^0$ has a quadratic Gr\"obner basis. 

(ii) $QA_n^0$ is Koszul. 

(iii) The Hilbert polynomial of $QA_n^0$ is $P_{{\mathfrak{qtr}}_n}(t)$. 
\end{corollary}

\begin{proof} (i) Pick any ordering of monomials with sum of labels monotonically 
nondecreasing. It is well known that the Orlik-Solomon algebra 
$OS_n$ has a quadratic Gr\"obner basis with respect to this ordering, compiled of all the relations 
(see \cite{Yu}): the initial monomials for this basis are products $b_{ip}b_{jp}$ with $p>i,j$. 
Putting this basis together with the elements 
$a_{ij}a_{jk}-a_{jk}a_{ki}$ for $k<i,j$, 
$a_{ij}b_{ij}$, $a_{ij}b_{jk}-a_{ik}b_{jk}$ for $k<j$, 
we get a quadratic Gr\"obner basis of $QA_n^0$. This implies (i).

(ii) follows from (i). 

(iii) This reduces to counting 2-step partitions $S$ with weights. 
We will adopt the following construction of such partitions: first we partition $[n]$ into 
$r$ nonempty subsets $S(i)$, $i=1,\dots ,r$, and then pick a set partition of $[r]$
into $l$ parts $T_1,\dots ,T_l$ to decide when we will have $S(i)=S_{pq}\subset S_p$;
namely, $S(i)\subset S_p$ if and only if $i\in T_p$.  

Let $s_p=|T_p|$, $p=1,\dots ,l$, and 
$d_i$ be the sizes of the parts $S(i)$.

Let $P(t)$ be the Hilbert polynomial of $QA_n^0$.
Let 
$$
F(t,u):=1+\sum_{n\ge 1}P_n(t)\frac{u^n}{n!}.
$$
We have 
\begin{multline}
F(t,u)=1+\sum_{r,\ell}\sum_{\substack{d_1,\dots,d_r>0\\d_1+\dots+d+r=n}}\frac{n!}{d_1!\cdots d_r!}\frac{(d_1-1)!\dots (d_r-1)!}{r!}\\
\times\sum_{\substack{s_1,\dots,s_\ell>0\\s_1+\dots +s_\ell=r}}\frac{r!}{s_1!\cdots s_\ell!}\frac{t^{n-\ell}}{\ell!}\frac{u^n}{n!}.
\end{multline}
Here $\frac{n!}{d_1!\cdots d_r!}$ is the number of ways to choose 
the parts $S(i)$ once the sizes of $S(i)$ have been fixed,
the factor $1/r!$ accounts for the fact that the order of 
the parts $S(i)$ does not matter, $(d_i-1)!$ are the sizes of the top components 
of the algebras $OS(S(i))$, 
$\frac{r!}{s_1!\dots s_\ell!}$ is the number of ways to choose 
the parts $T_p$ once their sizes have been fixed,
and $1/\ell!$ accounts for the fact that the order of 
the parts $T_p$ does not matter. 

Cancelling $n!,r!,(d_i-1)!$ and summing over $n,d_i$, we get 
$$
F(t,u)=1+\sum_{r,\ell}\sum_{\substack{s_1,\dots,s_\ell>0\\s_1+\dots+s_\ell=r}}
\frac{(-\log(1-tu))^r}{s_1!\dots s_\ell!}\frac{t^{-\ell}}{\ell!}.
$$
Now summing over $r>0,s_p$ we get 
$$
F(t,u)=1+\sum_{\ell}
\big((1-tu)^{-1}-1\big)^{\ell}\frac{t^{-\ell}}{\ell!}=e^{\frac{u}{1-tu}}.
$$
This completes the proof. 
\end{proof}

{\bf Remark.} Let $P(t,v)=\sum D_{pq}t^pv^q$, where 
$D_{pq}$ is the dimension of the space of elements of $QA_n^0$ of degree $p$, which have degree $n-q$ with 
respect to the variables $b_{ij}$. Set $F(t,u,v)=1+\sum_n P(t,v)u^n/n!$. 
Then the expression for $F(t,u,v)$ is obtained as above, except 
that we need to insert a factor $v^r$. 
This implies that 
$$
F(t,u,v)=\exp\frac{(1-tu)^{-v}-1}{t}
$$

\begin{proposition}\label{gr} 
(i) The natural map $QA_n^0\to {\rm gr}(QA_n)$ is an isomorphism. 

(ii) $QA_n$ is Koszul.  
\end{proposition}

\begin{proof}
(i) By Koszulity of $QA_n^0$, it is sufficient to check (see e.g.\ \cite{BG})
that the map is bijective in degrees $\le 3$, which is a direct computation.
Part (ii) follows from (i), since if ${\rm gr}(A)$ is Koszul, so is $A$. 
\end{proof}

Similarly to the previous section, 
Proposition \ref{gr} implies Theorem \ref{main} in the quasitriangular case.
Thus Theorem \ref{main} is proved. 

\subsection{Connection with the pure braid groups}

Let $\mathbf{PB}_n$ be the pure braid group on $n$ strands. 
Let ${\mathfrak{pb}}_n$ be the Lie algebra of its prounipotent
completion. According to the results of Kohno \cite{Ko},
this Lie algebra is isomorphic to its graded, and 
is generated by $t_{ij}=t_{ji}$, $i\ne j$, $i,j\in [n]$, 
with defining relations 
$$
[t_{ij},t_{ik}+t_{jk}]=0,
$$
and $[t_{ij},t_{kl}]=0$ if $i,j,k,l$ are distinct. 

We have a group homomorphism 
$\Psi_n: \mathbf{PB}_n\to {\mathbf{QTr}_n}$ defined by 
$$
T_{ij} \mapsto
R_{j-1,j}...R_{i+1,j} R_{ij}R_{ji} (R_{j-1,j}...R_{i+1,j})^{-1}
$$
where $T_{ij} = (\sigma_{j-1}... \sigma_{i+1})
\sigma_{i}^2 (\sigma_{j-1}... \sigma_{i+1})^{-1}$
are the Artin-Burau generators of $\mathbf{PB}_n$, and $\sigma_i$ are the Artin
generators of the full braid group.

The infinitesimal analog of $\Psi_n$ is the Lie algebra 
homomorphism $\psi_n: {\mathfrak{pb}}_n\to {\mathfrak{qtr}}_n$ 
defined by the formula $t_{ij} \mapsto r_{ij} +
r_{ji}$. 

It is clear that the kernel of the 
projection ${\mathbf {QTr}}_n\to {\mathbf {Tr}}_n$ is 
the normal closure of the image of $\Psi_n$. 
Similarly, the kernel of the 
projection ${\mathfrak{qtr}}_n\to {\mathfrak{tr}}_n$ is 
the normal closure of the image of $\psi_n$. 

\begin{proposition}\label{inje} 
The homomorphisms $\psi_n, \Psi_n$ are injective. 
\end{proposition}

The proof of Proposition \ref{inje} is given below in subsection 
\ref{proofinje}. 

{\bf Remark.} Here is another proof of the fact that $\psi_n$ is
injective. Let ${\mathfrak{qtr}}_n^0$ be the associated graded
of ${\mathfrak{qtr}}_n$ under the filtration defined by 
$\deg(t_{ij})=0$, $\deg(\rho_{ij})=1$. 
Proposition \ref{gr} implies that the defining relations for
this Lie algebra are
$$
[t_{ij},t_{ik}+t_{jk}]=0,\ [t_{ij},\rho_{ik}+\rho_{jk}]=0,
$$
$$
[\rho_{ij},\rho_{jk}]+
[\rho_{jk},\rho_{ki}]+[\rho_{ki},\rho_{ij}]=0,
$$
for distinct $i,j,k$, and
$$
[\rho_{ij},\rho_{kl}]=[\rho_{ij},t_{kl}]=[t_{ij},t_{kl}]=0
$$
if $i,j,k,l$ are distinct, and its universal enveloping algebra 
is the quadratic dual to $QA_n^0$. 
Let $\psi_n^0={\rm gr}\psi_n: {\mathfrak{pb}}_n\to
{\mathfrak{qtr}}_n^0$. Then it follows from the above relations
that $\psi_n^0$ (unlike $\psi_n$) is a {\it split} homomorphism,
hence it is injective. Thus $\psi_n$ is injective as well.


\newcommand{\univ}{\operatorname{univ}}
\newcommand{\cF}{\mathcal{F}}
\newcommand{\cA}{\mathcal{A}}
\newcommand{\SG}{\mathfrak{S}}
\newcommand{\kk}{{\bf k}}
\renewcommand{\ll}{{\bf l}}
\renewcommand{\wh}{\widehat}
\newcommand{\on}{\operatorname}
\newcommand{\ulI}{\underline{I}}
\newcommand{\ulJ}{\underline{J}}

\section{Proof of Proposition \ref{inje}}\label{proofinje}

In this section, we fix an integer $n\geq 1$. 
In \cite{En1}, the second author introduced 
$$
U^{\univ}_n = (U(\g)^{\otimes n})_{\univ} := 
\bigoplus_{N\geq 0} 
\big( (\cF\cA_N^{\otimes n})_{\sum_i \delta_i}
\otimes 
(\cF\cA_N^{\otimes n})_{\sum_i \delta_i}
\big)_{\SG_N}. 
$$ 
Here $\cF\cA_N$ is the free algebra generated by $x_1,\ldots,x_N$. 
It is graded by $\oplus_{i=1}^N \ZZ_{\geq 0}\delta_i$ 
(deg$(x_i) = \delta_i$), and equipped with the 
action of $\SG_N$ permuting the generators. The index $\sum_i\delta_i$
means the part of degree $\sum_i \delta_i$, and the index $\SG_N$
means the coinvariants of the diagonal action of $\SG_N$. 
We also defined an algebra structure on 
$U_n^{\univ}$. It has the following property: if 
$(A,r_A)$ is an algebra equipped with a solution $r_A = \sum_{i\in I} 
a(i) \otimes b(i) \in A^{\otimes 2}$ of the classical Yang-Baxter equation, 
then the linear map  $U_n^{\univ} \to A^{\otimes n}$ taking  
\begin{align*}
& \Big( x_1\cdots x_{k_1} \otimes x_{k_1+1}\cdots x_{k_2} 
\otimes \cdots \otimes x_{k_{n-1}+1}\cdots x_{N} \Big) 
\\ & \otimes 
\Big( x_{\sigma(1)}\cdots x_{\sigma(\ell_1)} \otimes x_{\sigma(\ell_1+1)}\cdots 
x_{\sigma(\ell_2)} \otimes \cdots \otimes 
x_{\sigma(\ell_{n-1}+1)}\cdots x_{\sigma(N)} \Big) 
\end{align*}
to  
\begin{align*}
& \sum_{i_1,\ldots,i_N\in I} a(i_1)\cdots a(i_{k_1})b(i_{\sigma(1)})
\cdots b(i_{\sigma(\ell_1)})
\otimes a(i_{k_1+1}) \cdots a(i_{k_2})  b(i_{\sigma(\ell_1 + 1)})
\cdots b(i_{\sigma(\ell_2)})
\\ & \otimes \cdots 
\otimes a(i_{k_{n-1}+1})\cdots a(i_N)
b(i_{\sigma(\ell_{n-1}+1)})
\cdots b(i_{\sigma(N)}) 
\end{align*}
is an algebra morphism. 

Let us set $r^{\univ}_{ij} = (1^{\otimes i-1} \otimes x_1 \otimes 1^{\otimes n-i}) 
\otimes (1^{\otimes j-1} \otimes x_1 \otimes 1^{\otimes n-j})
\in U_n^{\univ}$. 
Then we have a Lie algebra morphism $\zeta_n: \mathfrak{qtr}_n \to U^{\univ}_n$, 
$r_{ij}\mapsto r_{ij}^{\univ}$. 

By \cite{En2}, Section 1.13, the
composition $\zeta_n\circ \psi_n$ is injective.
Hence $\psi_n$ is injective.

Let $Z_n: {\mathbf{QTr}}_n\to (\wh U_n^{\rm univ})^\times_1$ be the
homomorphism sending $R_{ij}$ to $R_{ij}^{\rm univ}$. We denote by
$G(\QQ) = \operatorname{exp}(\operatorname{Lie}(G))$ the Malcev
$\QQ$-completion of a group $G$. As $(\wh U_n^{\rm univ})^\times_1$ is
a prounipotent $\QQ$-Lie group, $Z_n$ factors through $Z_n(\QQ):
\mathbf{QTr}_n(\QQ)\to (\wh U_n^{\rm univ})^\times_1$.  Also, let
$\Psi_n(\QQ): {\mathbf{PB}}_n(\QQ)\to {\mathbf {QTr}}_n(\QQ)$ be the
extension of $\Psi_n$ to Malcev $\QQ$-completions. It is easy to see
that the associated graded map of $Z_n(\QQ)\circ \Psi_n(\QQ):
\mathbf{PB}_n(\QQ)\to (\wh U_n^{\rm univ})^\times_1$ is
$\zeta_n\circ\psi_n$. Hence $Z_n(\QQ)\circ \Psi_n(\QQ)$ is injective,
and therefore $\Psi_n(\QQ)$ is injective. But the group
$\mathbf{PB}_n$ is an iterated cross product of free groups, which
implies that the natural map $\mathbf{PB}_n\to \mathbf{PB}_n(\QQ)$ is
injective.  Thus $\Psi_n$ is injective, as desired.

\begin{question}
Is the map $\zeta_n$ injective? 
\end{question}

\section{Classifying spaces for the groups ${\mathbf {Tr}}_n$,
${\mathbf {QTr}}_n$.}

\subsection{The Permutohedron}
Let $P_n$ be the convex hull of $\sym_n\cdot(n,n-1,\dots,1)$ in the
affine hyperplane $\mathbb A_n$ defined by the equation
$$\sum_{i=1}^n x_i=1+\dots+n=n(n+1)/2$$
in $\mathbb R^n$. This is a polyhedron, containing the points
$(x_1,\dots,x_n)$ such that, for every $S\subset[n]$, we have
$\sum_{s\in S}x_s\in[1+\dots+|S|,(n-|S|+1)+\dots+n]$.

For $S$ a finite set, we write $P_S$ for the permutohedron $P_{|S|}$
constructed in $\mathbb R^S$.

The faces of $P_n$ can be determined as follows. The
$(n-r)$-dimensional faces are in bijection with the ordered partitions
$[n]=S_1\sqcup\dots\sqcup S_r$; the face corresponding to such a
choice is the set of the points $(x_1,\dots,x_n)$ satisfying
$$\text{for all }i\in\{1,\dots,r\}:\quad\sum_{s\in S_i}x_s=(|S_1|+\dots+|S_{i-1}|+1)+\dots+(|S_1|+\dots+|S_i|).$$
This face is therefore the Cartesian product $P_{S_1}\times\dots\times
P_{S_r}$, with the coordinates of $P_{S_i}$ shifted
$|S_1|+\dots+|S_{i-1}|$ away from the origin. For example, the vertex
with coordinates $(\pi(1),\dots,\pi(n))$ corresponds to the partition
$[n]=\{\pi^{-1}(1)\}\sqcup\dots\sqcup\{\pi^{-1}(n)\}$.  Geometric
inclusion of faces corresponds combinatorially to ordered refinement
of partitions; namely, $S_1\sqcup\dots\sqcup S_r$ is a face of
$T_1\sqcup\dots\sqcup T_s$ if there is an order-preserving surjection
$f:[r]\to[s]$ with $T_i=\bigcup_{j\in f^{-1}(i)}S_j$. The set of faces
is partially ordered by this relation; partitions into singletons are
atoms, and the one-part partition $[n]$ is the maximal element.

\subsection{The classifying space $C_n$}
For every $r$ there is a natural action of the symmetric group
$\sym_r$ on the disjoint union of the $(n-r)$-dimensional faces of
$P_n$. In the combinatorial model, it is the natural permutation
action
$$\pi(S_1\sqcup\dots\sqcup S_r)=S_{\pi^{-1}(1)}\sqcup\dots\sqcup S_{\pi^{-1}(r)},$$
which moves $S_i$ from position $i$ to position $\pi(i)$. In the
geometric model, the action is given by piecewise translations: in the
face associated with $S_1\sqcup\dots\sqcup S_r$,
$$\pi(x_s)=x_s-\left(|S_1|+\dots+|S_{i-1}|\right)+\left(|S_{\pi^{-1}(1)}|+\dots+|S_{\pi^{-1}(\pi(i)-1)}|\right).$$
These actions fit together, in the sense that if $S_1\sqcup\dots\sqcup
S_r$ is a face of $T_1\sqcup\dots\sqcup T_s$ via the surjection $f$ as
above, then $f$ interlaces the $\sym_r$- and $\sym_s$-actions. We let
$C_n$ be the quotient of $P_n$ by these actions.

\begin{theorem}\label{thmCS}
  $C_n$ is a classifying space for the group ${\mathbf {Tr}}_n$.
\end{theorem}
\begin{proof}
  This amounts to showing that $C_n$ has a contractible universal
  cover, and has fundamental group ${\mathbf {Tr}}_n$.

  For the first assertion, we show that $C_n$ is locally a
  non-positively curved space, whence~\cite[Corollary~II.1.5]{BH}
  yields that the universal cover of $C_n$ is contractible. Since
  $C_n$ is a quotient of Euclidean space, it suffices to show
  by~\cite[Theorem~II.5.2]{BH} that the link $Lk(*,C_n)$ of the vertex
  $*\in C_n$ has curvature $\le1$; applying
  iteratively~\cite[Theorem~5.4]{BH}, this amounts to showing that for
  every face $F$ of $C_n$ the link $Lk(F,C_n)$ contains no
  isometrically embedded circles of length $<2\pi$.

  Recall that the faces of $P_n$ are indexed by ordered set-partitions
  of $[n]=\{1,2,\dots,n\}$, and that faces of $C_n$ are indexed by
  unordered set-partitions of $[n]$.

  Let therefore $\mathcal S=\{S_1,S_2,\dots,S_r\}$ be an unordered set
  partition, with $S_i \subset [n]$, and let $F$ denote the
  corresponding face of $C_n$. Note that $F$ is of dimension $n-r$.
  The link $L=Lk(F_S,C_n)$ is a spherical simplicial complex of
  dimension $r-1$, combinatorially isomorphic to $Lk(*,C_r)$, which
  can be described as follows:

  The vertex set $V$ of $L$ consists of ordered pairs of distinct
  elements of $S$; the pair $(S_i,S_j)$ represents the partition
  $\{S_1,\dots,S_i\cup S_j,\dots,S_r\}$. A subset $\Theta$ of $V$
  spans a simplex if and only if there exists a permutation
  $(T_1,\dots,T_r)$ of $\mathcal S$ such that
  $\Theta\subset\{(t_1,t_2),(t_2,t_3),\dots,(t_{r-1},t_r)\}$.

  An isometrically embedded circle in $L$ is necessarily a subset of
  the $1$-skeleton of $L$, so it is important to compute the lengths
  of the edges of $L$. These lengths can be described geometrically as
  dihedral angles between certain $(n-r+1)$-faces incident to $F$. By
  considering normal vectors to these $(n-r+1)$-faces inside the
  $(n-r+2)$-face in $\mathbb A_n$ containing them, these angles can be
  computed explicitly.
  
  Let $(S_i,S_j)$ and $(S_k,S_l)$ be two vertices of $L$. For
  $S\subset[n]$, let $e_S$ denote the vector having $1$'s in
  co\"ordinates $\in S$ and $0$ elsewhere. Two cases can occur:
  \begin{itemize}
  \item $\boldsymbol{|\{i,j,k,l\}|=4}$: normal vectors to $(S_i,S_j)$ and
    $(S_k,S_l)$ may be chosen respectively as
    $|S_j|e_{S_i}-|S_i|e_{S_j}$ and $|S_k|e_{S_l}-|S_l|e_{S_k}$, so
    that the arclength between these two vertices is $\pi/2$.
  \item $\boldsymbol{j=k,|\{i,j,l\}|=3}$: normal vectors may be chosen
    respectively as
    $v_{ij}=(|S_k|+|S_l|)e_{S_i}-|S_i|e_{S_k}-|S_i|e_{S_l}$ and
    $v_{kl}=(|S_i|+|S_j|)e_{S_l}-|S_l|e_{S_i}-|S_l|e_{S_j}$, and the
    arclength between these vertices is the angle between these
    vectors. The exact value will not be needed, but since $\langle
    v_{ij},v_{kl}\rangle<0$ the angle is strictly $>\pi/2$.
  \end{itemize}

  Consider now an isometrically embedded circle in $L$; we wish to
  show that its length it at least $2\pi$. Since all edge lengths in
  the $1$-skeleton of $L$ have length at least $\pi/2$, it suffices to
  consider triangles in $L$. The only geodesic triangles are those
  which do not bound a $2$-face. There is one such triangle for every
  cyclically ordered triple of elements of $S$; the vertices of this
  triangle are of the form $(S_i,S_j)$, $(S_j,S_k)$ and $(S_k,S_i)$.
  Now the arclength between these vertices add up to $2\pi$, because
  their respective vectors $v_{ij}$, $v_{jk}$ and $v_{ki}$ are
  coplanar: $v_{ij}+v_{jk}+v_{ki}=0$.

  We next compute the fundamental group of $C_n$. It is generated by
  simple loops in the $1$-skeleton of $C_n$, and has relations given
  by the $2$-cells. A simple loop in the $1$-skeleton is of the form
  $\{a,b\}$, neglecting singletons. We identify it with the generator
  $R_{ab}$, with the ordering $a<b$. The other generators $R_{ba}$ are
  redundant, because of the relation $R_{ab}R_{ba}=1$.

  A typical $2$-cell is either of the form $\{a,b\}\sqcup\{c,d\}$, in
  which case it gives the relation $R_{ab}R_{cd} = R_{cd}R_{ab}$, 
  or of the form
  $\{a,b,c\}$, in which case it gives the relation
  $R_{ab}R_{ac}R_{bc}=R_{bc}R_{ac}R_{ab}$, again assuming the ordering
  $a<b<c$.

  The other relations of $\mathbf {Tr}_n$, namely those of the form
  $R_{ab}R_{ac}R_{bc}=R_{bc}R_{ac}R_{ab}$ with $a,b,c$ not in the
  order $a<b<c$, are cyclic permutations of the one in the standard
  ordering.
\end{proof}

\begin{remark}
  The quotient $C_n$ of $P_n$ may be constructed in two steps. First,
  let $T_n$ be the space obtained from $P_n$ by identify opposing
  faces $S_1\sqcup S_2$ and $S_2\sqcup S_1$. In terms of partitions,
  this corresponds to identifying an ordered partition with all of its
  cyclic permutations. There is a lattice of translations of $\mathbb
  A_n$, isomorphic to $\ZZ^{n-1}$, and spanned by all vectors
  $(-1,\dots,-1,n-1,-1,\dots,-1)$. It is easy to see that $P_n$ is a
  fundamental domain for this lattice. Indeed two translates with
  non-trivial intersection are of the form $P_n+v$ and
  $P_n+v+|S|\sum_{s\not\in S}x_s-(n-|S|)\sum_{s\in S}x_s$, for some
  $S\subset[n]$, and these two translates intersect in the
  $(n-2)$-dimensional face $S\sqcup([n]\setminus S)$.

  There are natural copies of $P_{i_1}\times\dots\times P_{i_r}$ with
  $i_1+\dots+i_r=n$ in $T_n$, for example a copy of $P_{n-1}\times
  P_1\cong P_{n-1}$ spanned by all $S_1\sqcup\dots\sqcup
  S_r\sqcup\{n\}$ for any face $S_1\sqcup\dots\sqcup S_r$ of
  $P_{n-1}$. One may construct $C_n$ by quotienting in $T_n$ each of
  these $P_{i_1}\times\dots\times P_{i_r}$ into
  $C_{i_1}\times\dots\times C_{i_r}$.

  These considerations are sufficient to describe the spaces $C_n$ for
  small $n$: $C_1$ is a point; $C_2$ is a circle, the quotient of the
  line segment $[(1,2),(2,1)]$ by identification of its endpoints; and
  $C_3$ is the quotient of a $2$-torus obtained by gluing two distinct
  points together. It is homotopic to the connected sum of a $2$-torus
  and a circle.
\end{remark}

\subsection{The classifying space $QC_n$.}

Consider now the space $QC_n$ constructed as follows. On the disjoint
union $C_n\times\sym_n$, identify all faces $(S_1\sqcup\dots\sqcup
S_r,\sigma)$ and $(S_1\sqcup\dots\sqcup S_r,\tau)$ precisely when
\[\text{ for all }i\in\{1,\dots,r\}:\,\text{ for all }x,y\in
S_r:\,\sigma(x)<\sigma(y)\Leftrightarrow \tau(x)<\tau(y).\]

\begin{theorem}
  $QC_n$ is a classifying space for the group $\mathbf {QTr}_n$.
\end{theorem}
\begin{proof}
  The proof that $QC_n$ is locally a non-positively curved space is
  similar to the proof for $C_n$: the links in $QC_n$ are covers of
  links in $C_n$. We omit the details.

  A simple loop in the $1$-skeleton of $QC_n$ is of the form
  $(\{a,b\},\sigma)$, neglecting singletons, where $\sigma$ specifies
  an orientation of $\{a,b\}$. We identify it with the generator
  $R_{ab}$, with the ordering specified by $\sigma(a)<\sigma(b)$.

  A typical $2$-cell is either of the form
  $(\{a,b\}\sqcup\{c,d\},\sigma)$, in which case it gives the relation
  $R_{ab}R_{cd}=R_{cd}R_{ab}$, in the order specified by $\sigma(a)<\sigma(b)$
  and $\sigma(c)<\sigma(d)$, or is of the form $(\{a,b,c\},\sigma)$,
  in which case it gives the relation
  $R_{ab}R_{ac}R_{bc}=R_{bc}R_{ac}R_{ab}$, again in the order
  $\sigma(a)<\sigma(b)<\sigma(c)$.
\end{proof}

Now the fact that $\mathbf{Tr}_n$ is a split quotient of
$\mathbf{QTr}_n$ acquires a geometric interpretation.
Namely, $C_n$ is both a quotient of $QC_n$, by
projecting on the first coordinate, and a subspace of $QC_n$,
embedded as $C_n\times\{1\}$.

\subsection{Homology of $C_n$ and $QC_n$}
The permutohedron $P_n$ is homeomorphic to a ball, and therefore has
non-trivial homology only in dimension $0$. Consider the
chain complex spanned by all faces of the polyhedron. In the 
combinatorial model, the boundary operator is given by
\begin{align}
  \partial(S_1\sqcup\dots\sqcup S_r)&=\sum_{i=1}^r(-)^i
  S_1\sqcup\dots\sqcup\partial(S_i)\sqcup\dots\sqcup S_r,\label{eq:bdry1}\\
  \partial(S)&=\sum_{S=S'\sqcup S''}(-)^? S'\sqcup S'',\notag
\end{align}
for appropriate signs in $\partial(S)$. We fix an orientation for the
top-dimensional face in $P_n$. Each face $S_1\sqcup\dots\sqcup S_r$ of
$P_n$ is then a Cartesian product of translates of $P_{S_i}$, and we
give that face the product orientation. The exact signs in
$\partial(S)$ are not important, but we note that $S'\sqcup S''$ and
$S''\sqcup S'$ are translates of each other, and have opposite signs
in $\partial(S)$.

Consider now the quotient $C_n$. It admits a natural chain complex,
where the $(n-r)$-dimensional complex is spanned by $\sym_r$-orbits of
$(n-r)$-dimensional faces, i.e.\ by unordered partitions of $[n]$ in
$r$ parts. The boundary operator on that complex is induced from the
boundary on $P_n$, since by the choice above all faces in the same
$\sym_r$-orbit have the same orientation.

\begin{proposition}
(i)  The boundary map on $P_n$ induces the zero map on the complex of
  $C_n$; in other words, the complex for $C_n$ is minimal.

(ii) The homology group $H_r(\mathbf{Tr}_n;\ZZ)$ is free
  of rank the number $\{\begin{smallmatrix}n\\ n-r\end{smallmatrix}\}$
  of unordered partitions of $n$ in $n-r$ parts\footnote{a.k.a. the
  ``Stirling numbers of the second kind''}. The Hilbert-Poincar\'e
  series of $H_*(\mathbf{Tr}_n;\ZZ)$ is the polynomial
  $P_{{\mathfrak{tr}}_n}(t)$ of Theorem~\ref{main}.
\end{proposition}
\begin{proof}
  The boundary map on any face of $C_n$ is calculated from the
  boundary on top-dimensional faces, by formula~\eqref{eq:bdry1}. Now
  the top-dimensional face of $P_n$ is mapped bijectively to $C_n$,
  while the faces of dimension one less are identified by pairs in
  $C_n$. Furthermore, the boundary operator assigns sign $+1$ to one
  and sign $-1$ to the other. Therefore all boundary operators are
  trivial on top-dimensional faces, and by extension on all faces.

  The consequence for homology of $\mathbf{Tr}_n$ follows immediately,
  since it is by definition the homology of a classifying space.
\end{proof}

The same arguments hold for $QC_n$, namely:
\begin{proposition}
(i)  The boundary map on $P_n$ induces a trivial map on the complex of
  $QC_n$; in other words, the complex for $QC_n$ is minimal.

(ii) The homology group $H_r(\mathbf{QTr}_n;\ZZ)$ is free
  of rank the number of unordered partitions of $n$ in $r$ ordered
  parts.  The Hilbert-Poincar\'e series of $H_*(\mathbf{QTr}_n;\ZZ)$
  is the polynomial $P_{{\mathfrak{qtr}}_n}(t)$ of
  Theorem~\ref{main}.\end{proposition}
\begin{proof}
  Since $QC_n$ is a quotient of a disjoint union of copies of $C_n$,
  the induced boundary map is \emph{a fortiori} trivial.

  The $r$-faces of $QC_n$ are given by the numbers $A_{n,n-r}$ of
  partitions of $[n]$ in $n-r$ parts, with an ordering on each of the
  parts. These numbers obviously satisfy the recursion
  \[A_{n,p}=A_{n-1,p-1}+(n+p-1)A_{n-1,p},\]
  since given a partition of $[n]$ in $p$ parts one may remove $n$
  from the partition and obtain either a partition of $[n-1]$ in $p-1$
  parts, or a partition of $[n-1]$ in $p$ parts, where the number $n$
  appeared in any of $n+p-1$ positions. This recursion is also
  satisfied by the coefficients $\binom{n-1}{p}\frac{n!}{(n-p)!}$
  appearing in Theorem~\ref{main}, finishing the proof.
\end{proof}

\subsection{The cohomology rings of the groups ${\mathbf {Tr}}_n$,
${\mathbf {QTr}}_n$.}

We now relate the cohomology of the Lie algebras 
${\mathfrak{(q)tr}}_n$ and the groups ${\mathbf{(Q)Tr}}_n$.

We have natural homomorphisms $\xi_n: A_n\to H^*({\mathbf
  {Tr}}_n,\QQ)$ and $\eta_n: QA_n\to H^*({\mathbf {QTr}}_n,\QQ)$,
defined as follows. The generators $a_{ij}$ of $H^*({\mathbf
  {Tr}}_n,\QQ)$ and $a_{ij},b_{ij}$ of $H^*({\mathbf {QTr}}_n,\QQ)$
are obtained by pulling back standard generators via the projections
${\mathbf {Tr}}_n\to {\mathbf {Tr}}_2$, ${\mathbf {QTr}}_n\to {\mathbf
  {QTr}}_2$ (note that the groups ${\mathbf {Tr}}_2$, ${\mathbf
  {QTr}}_2$ are free in one, respectively two generators). The fact
that these generators satisfy the relations of $A_n, QA_n$
respectively follows from consideration of the projections ${\mathbf
  {Tr}}_n\to {\mathbf {Tr}}_3$, ${\mathbf {QTr}}_n\to {\mathbf
  {QTr}}_3$, and the structure of the cohomology rings of ${\mathbf
  {Tr}}_3$, ${\mathbf {QTr}}_3$, which is easy to determine by looking
at the 2-dimensional complexes $C_3,QC_3$. This defines the
homomorphisms $\xi_n,\eta_n$.

\begin{conjecture}\label{cohring} (stated as Theorem 8.5 in the 
published version; now a theorem of P. Lee, \cite{L})
$\xi_n$ is an isomorphism. 
\end{conjecture}

Let us show that Conjecture \ref{cohring} follows from Conjecture \ref{isomo}.
Indeed, by the above results, $\xi_n$ 
is a map between spaces of the same dimension. 
Thus it suffices to show that $\xi_n$ is injective. 

Let $\g_n=\operatorname{Lie}{\mathbf{Tr}}_n$. We have a natural map
$\theta_n: H^*({\mathbf{Tr}}_n,\QQ)\to H^*_{cts}(\g_n,\QQ)$ (the
subscript cts means continuous cohomology).  By Conjecture \ref{isomo},
we have ${\rm gr}\g_n=\g_n^0={\mathfrak{tr}}_n$.  Thus, since the
algebra $U(\g_n^0)$ is Koszul, all the differentials of the spectral
sequence computing the cohomology of $H^*_{cts}(\g_n,\QQ)$ starting
from $H^*(\g_n^0,\QQ)$ are zero.  Thus the injectivity of $\xi_n$ will
follow from the injectivity of the natural map $\xi_n^0: A_n\to
H^*(\g_n^0,\QQ)$ (as $\theta_n\circ \xi_n$ is a deformation of
$\xi_n^0$).  But we know (by the results of Sections 3,4) that
$\xi_n^0$ is an isomorphism, as desired. 

The same argument shows that Conjecture
\ref{isomo1} implies the following 

\begin{conjecture}\label{cohring1} (now a theorem of P. Lee, \cite{L})
$\eta_n$ is an isomorphism. 
\end{conjecture}

\subsection*{Acknowledgments} The authors are very grateful to:
Richard Stanley for useful discussions; Andr\'e Henriques for contributing a
proof of Theorem~\ref{thmCS}; Leonid Bokut and Peter Lee for pointing out 
errors in the previous version and for giving useful references.  P.E.\
and B.E.\ thank the mathematics department of ETH (Z\"urich) for
hospitality. The work of P.E. was partially supported by the NSF grant
DMS-0504847 and the CRDF grant RM1-2545-MO-03. E.R.\ was supported in
part by NSF Grant No.  DMS-0401387. Throughout the work, we used the
``Magma'' program for algebraic computations, \cite{Ma}.

\end{document}